\documentclass{amsart}

\title[Fatou-Bieberbach Domains]{Fatou-Bieberbach Domains}
\author{Erlend Forn\ae ss Wold}

\date{October 5, 2004}

\usepackage{amscd}

\newtheorem{theorem}{Theorem}
\newtheorem{lemma}{Lemma}
\newtheorem{proposition}{Proposition}

\theoremstyle{definition}
\newtheorem{defin}{Definition}

\theoremstyle{remark}

\newcommand{\NN}{\mathbb{N}}
\newcommand{\RR}{\mathbb{R}}
\newcommand{\CC}{\mathbb{C}}

\def\C{{\bf C}}

\def\e{{\epsilon}}
\def\d{{\delta}}
\def\l{{\lambda}}
\def\r{{\rho}}

\def\cO{{\mathcal O}}

\begin{document}

\begin{abstract}
We show that for any $m\in\NN\cup\{\infty\}$  there exist $m$
disjoint FB domains whose union is dense in $\CC^k$.  In fact we
show that any point not in the union is a boundary point for all the
domains. We construct FB domains that contains arbitrary countable
collections of subvarieties of $\CC^k$, and we construct FB domains
that intersect elements of countable collections of affine subspaces
of $\CC^k$  in connected proper subsets. Moreover, we show that any
Runge FB domain is the attracting basin for a sequence of
automorphisms of $\CC^k$, although not necessarily if you only allow
iteration of  one automorphism.  We also show that an increasing
sequence of Runge $\CC^k$'s is a $\CC^k$.
\end{abstract}

\maketitle

\section{Introduction}
This paper is inspired by the paper [10], and is organized as
follows: We start by giving some definitions in Section 2. In
Section 3 we will prove a generalization of Theorem 9.1 in [10]
that, together with results on approximations of biholomorphisms by
automorphisms of $\CC^k$, due to Andersen, Lempert, Forstneric, and
Rosay, will be a very effective tool for constructing various
Fatou-Bieberbach Domains. In Section 4 we develop some results
regarding polynomial convexity needed in Section 6.  In Section 5 we
prove the following theorem: \

\begin{theorem} For any $m\in\NN\cup\{\infty\}$ there exists a union
of Fatou-Bieberbach Domains
 $\Omega=\cup_{j=1}^m\Omega_j$  that satisfies the following:

\

(i) $\Omega_k\cap\Omega_l=\emptyset$  for all $k\neq l$, \

(ii) For any $q\in\CC^k\setminus\Omega$  we have that
$q\in\partial\Omega_i$  for all i.
\end{theorem}

\

This answers a question posed by Rosay and Rudin in [10]. Rosay and
Rudin also posed a couple of questions regarding intersections
between Fatou-Bieberbach domains an complex lines in $\CC^2$. One
question is the following: Can the intersection between a
Fatou-Bieberbach domain and a complex line in $\CC^2$  be connected?
Globevnik has given a positive answer to this question by
constructing a Fatou-Bieberbach domain intersecting the $z$-plane in
approximately a disc [8].  Another question is the following: How
many complex lines can a Fatou-Bieberbach domain in $\CC^2$ contain?
Buzzard and Forn\ae ss have shown that a Fatou-Bieberbach domain in
$\CC^2$  can contain any finite number of complex lines [2].  In
Section 6 we prove the following two theorems:\

\label{affine}\begin{theorem} Let $\{L_j\}_{j\in\NN}$  be a
collection of affine subspaces of $\CC^k$.  Then there exists a
Fatou-Bieberbach Domain $\Omega$  such that $\Omega\cap L_i$  is
connected, and such that $L_i\setminus\Omega\neq\emptyset$  for all
$i\in\NN$.
\end{theorem}
\label{subvariety}\begin{theorem} Let $\{V_j\}_{j\in\NN}$  be a
collection of closed subvarieties of $\CC^k$. Then there exists a
Fatou-Bieberbach Domain $\Omega$ such that $\cup_{j=1}^\infty
V_j\subset\Omega$.
\end{theorem}
Notice that $\{V_j\}$  could be dense in $\CC^k$.  Theorem 3
generalizes a result in [10] stating that a Fatou-Bieberbach domain
can contain any countable set of points. \

In Section 7 we give an example of a Fatou Bieberbach domain that is
not the basin of attraction of an automorphism of $\CC^k$, but on
the other hand we prove that any Runge Fatou-Bieberbach domain is
the basin of attraction for a sequence of automorphism of $\CC^k$.
Lastly, we show that an increasing union of Runge domains that are
biholomorphic to $\CC^k$  is again a $\CC^k$.  This last result
gives a partial answer to a question posed in [7, p.4].

\

Apart from being interesting in their own right, constructions of
Fatou-Bieberbach domains with special properties can have useful
applications in other areas.  See [11] for an application to proper
holomorphic embeddings. \

The author would like to thank the referee for useful comments and
suggestions. \

\section{Definitions and Notation}

Throughout the article we will use the following notation: If
$p\in\CC^k$  and if $\e>0$, we let $B_\e(p)$  denote open ball
centered at $p$  with radius $\e$.  If $k=1$  this set is denoted
$\triangle_\e(p)$.  If $p$  is the origin, these sets will be
denoted $B_\e$  and $\triangle_\e$, and if in addition $\e=1$, these
sets are denoted $B$ and $\triangle$ respectively.  If nothing else
is stated, k is always assumed to be larger than or equal to two. \

\begin{defin}
Let $Aut_p(\CC^k)$  denote the group of holomorphic automorphisms of
$\CC^k$  fixing the point $p\in\CC^k$.  If all the eigenvalues
$\l_i$ of $\mathrm{d}F(p)$  satisfy $|\l_i|<1$  we say $F$  is
attracting at $p$.
\end{defin}

\

\begin{defin}
Let $\{F_j\}\subset Aut_p(\CC^k)$.  We write $F(i,j)=F_j\circ
...\circ F_i$  for $i\leq j$, and if $i=1$  we write $F(j)$  for
short.  If $i>j$  we will let $F(i,j)$  be the identity.  We define
the basin of attraction of a point $p\in\CC^k$
 by
$$
\Omega_{\{F_j\}}^p=\{z\in\CC^k;\lim_{j\rightarrow\infty}
F(j)(z)=p\}.
$$
\end{defin}

\

\section{Fatou-Bieberbach Domains As Sequence Attracting Basins}

It was proven in [10] that $\Omega_{\{F\}}^p$  is a Fatou-Bieberbach
domain for all attracting $F\in Aut_p(\CC^k)$.  A weaker result was
also proven, one in which was put an additional condition on the
eigenvalues of $\mathrm{d}F(p)$. To help us construct our domains,
we will generalize this weaker result, and Theorem 4 should be
compared with Theorem 9.1 in [10].

\begin{lemma}
Let $\Gamma(s,r,\r)$  be the family of biholomorphic maps $F\colon
B_\r\rightarrow\CC^k$  that fixes the origin and satisfies
$s\|z\|\leq\|F(z)\|\leq r\|z\|$  for all $z\in B_\r$, where
$s,r,\r\in\RR^+$.  Let $A_F=dF(0)$.  There exists a $C>0$ such that
we for all $F\in\Gamma$  have that
$$
\|A_F^{-1}F(z) - z\|\leq C\|z\|^2
$$
for all $z\in B_\r$.
\end{lemma}

\begin{proof}
For any $F\in\Gamma(s,r,\r)$, the map $G=A_F^{-1}F - I$ has no terms
of order less than two. This means that there exists a $C$ for this
particular map. But as we can choose such a $C$ depending only on
the supremum of $G$ over $B_\r$, it follows from the boundedness of
$\Gamma(s,r,\r)$  on $B_\r$  that one single $C$ must work for all
maps.\
\end{proof}
\label{square}\begin{theorem}\ Let $0<s<r<1$  such that $r^2<s$, let
$\d>0$, and let $\{F_j\}\subset Aut_p(\CC^k)$ such that
$s\|z-p\|\leq\|F_j(z)-p\|\leq r\|z-p\|$ for all $z\in B_\d(p)$ and
all $j\in\NN$. Then there exists a biholomorphic map
$$
\Phi\colon\Omega^p_{\{F_j\}}\rightarrow\Phi(\Omega^p_{\{F_j\}})=\CC^k.
$$
\end{theorem}
\begin{proof}
We may assume that $p=0$  and that $\d<1$, and we write
$\Omega=\Omega^0_{\{F_j\}}$.  Let $A_j=\mathrm{d}F_j(0)$.  We will
prove that the sequence of automorphisms $\Phi_j$ defined by
$$
\Phi_j=A(j)^{-1}F(j)
$$ converges to the desired map
uniformly on compacts in $\Omega$.  For all $z\in B_\d$  and all
$j\in\NN$ we have that $\|F(j)(z)\|\leq r^j\|z\|<r^j$, so the
sequence $\{F(j)\}$  is uniformly attracting on $B_\d$.   It follows
that
$$
\Omega=\cup_{j=0}^\infty F(j)^{-1}(B_\d),
$$
from which it follows that $\Omega$  is an open and connected subset
of $\CC^k$.
\

To prove convergence of the sequence $\Phi_j$  on compacts in
$\Omega$ it is enough to prove convergence of $\Phi_j$  on $B_\d$.
For if $K\subset\Omega$  is a compact set it is clear that for a
$k\in\NN$ we have that $F(k)(K)\subset B_\d$, and for any $z\in K$,
the limit
$$
\Phi(z)=\lim_{j\rightarrow\infty}\Phi_j(z)
$$
can be written as
$$
\Phi(z)=\lim_{j\rightarrow\infty} A(k)^{-1}\Phi^k_j
F(k)(z)=A(k)^{-1}\Phi^k F(k)(z),
$$
where $\Phi^k_j$  is the composition
$$
\Phi^k_j=A(k+1,j)^{-1} F(k+1,j),
$$
and $\Phi^k$  is the limit
$$
\Phi^k=\lim_{j\rightarrow\infty}\Phi^k_j.
$$
Remember that $\{F_j\}$ is an arbitrary sequence of maps in
$\Gamma(s,r,\d)$.  Lemma 1 gives us the following estimate:
\begin{align*}
\|\Phi_{j+1}(z)-\Phi_j(z)\|&=\|A(j)^{-1}(A_{j+1}^{-1}
F_{j+1}(F(j)(z)) - F(j)(z))\|  \\
&\leq s^{-j}C\|F(j)(z)\|^2\leq C\d^2(r^2/s)^j
\end{align*}
for all $z\in B_\d, j\in\NN$.  Since $\sum_{j=0}^\infty (r^2/s)^j$
exists, this shows that the sequence is uniformly convergent, i.e.
$\Phi$  is a holomorphic map from $\Omega$  into $\CC^k$.  Since we
have that $\mathrm{d}\Phi_j(0)=I$  for all $j\in\NN$, we have that
$\mathrm{d}\Phi(0)=I$, and it follows from a standard result that
the limit map has to be one to one onto its image.  It remains to
show that the image is the whole of $\CC^k$.\

Notice that there exists an $R\in\RR^+$  such that
$$
(a) \ \Phi^k(B_\d)\subset B_R,
$$
for all $k\in\NN$. \

We claim that there exists an $\e>0$  such that for any $k\in\NN$ we
have that
$$
(b) \ B_\e\subset\Phi^k(B_\d).
$$
For if not there is a sequence of positive numbers $\e_i\searrow 0$
with a corresponding sequence $\{k_i\}$  such that
$$
(c) \ B_{\e_i}\setminus\Phi^{k_i}(B_\d)\neq\emptyset
$$
for all $i\in\NN$.  By $(a)$  we have that $\{\Phi^{k_i}\}$ is a
normal family on $B_\d$  so we may assume that
$$
\lim_{i\rightarrow\infty}\Phi^{k_i}=\tilde\Phi,
$$
where $\tilde\Phi$  is a biholomorphic map on $B_\d$  fixing the
origin. This leads to a contradiction as $\tilde\Phi(B_\d)$  clearly
would have to contain a ball of some positive radius, contradicting
(c). \

For an arbitrary $M\in\RR^+$  there exists a $k\in\NN$  such that
$B_M\subset A(k)^{-1}(B_\e)$.  It follows from $(b)$  that
$$
B_M\subset A(k)^{-1}(B_\e)\subset A(k)^{-1}(\Phi^k(B_\d))\subset
A(k)^{-1}\Phi^k F(k)(\Omega)=\Phi(\Omega),
$$
which means that $\CC^k=\Phi(\Omega)$; thus the proof is finished.
\end{proof}

\

In Section 7 we will prove that even though quite restrictive, the
above theorem is (theoretically) sufficient for constructing all
Runge Fatou-Bieberbach domains.  A non-Runge Fatou-Bieberbach domain
would obviously not be a basin of attraction.\

It is an open question whether the theorem would still hold if we
drop the condition $r^2<s$. It is however clear that one in general
needs some upper and lower bound for the family of automorphisms.
Without an upper bound one could choose a sequence of linear maps
approaching the identity so fast that the basin would simply be the
origin.  A more interesting example can be found in [6].  In this
paper, Forn\ae ss constructs sequences of automorphisms where the
lower bound decreases fast to zero.  The resulting basins are
increasing unions of holomorphic balls, and the infinitesimal
Kobayashi metric vanishes identically on the basins. Nevertheless
they fail to be biholomorphic to $\CC^k$ due to the fact that they
carry nonconstant bounded plurisubharmonic functions.

\

\section{Polynomial Convexity}

For our constructions in connection with subvarieties of $\CC^k$,
we will need some results concerning polynomial convexity.  These
results are however not needed in Section 5.\

The holomorphically convex hull of a compact set $K\subset U$,
with respect to the set $U$, is defined as
$$
\widehat K_{\cO(U)}=\{z\in U;\|f(z)\|\leq\|f\|_K,\forall f\in
\cO(U)\}
$$
If $U=\CC^k$  we suppress the subscript, and write $\widehat K$
instead of $\widehat K_{\cO(\CC^k)}$.  If for a compact set $K$ we
have that $K=\widehat K$, we say that $K$  is polynomially convex. \

\label{variety}\begin{lemma} Let $K\subset\CC^k$  be polynomially
convex, let $V\subset\CC^k$  be a closed subvariety, and let
$K'\subset V$ be compact such that $K\cap V\subset K'$. Then we have
that
$$
\widehat{K\cup K'}=K\cup\widehat K'_{\cO(V)}=K\cup\widehat
K'
$$
\end{lemma}
\begin{proof}
Write $C=K\cup K'$.  Let $p\in\CC^k\setminus(C\cup V)$. There is an
$f\in\cO(\CC^k)$  such that $\|f\|_K<1$  and such that $\|f(p)\|>1$,
and there is an $h\in\cO(\CC^k)$  such that $h\mid_V\equiv 0$  and
such that $h(p)\neq 0$.  So if we define $g_k(z)=h(z)\cdot f(z)^k$,
we have that $\|g_k(p)\|>\|g\|_{C}$ for a large enough $k$.  It
follows from the Local Maximum Modulus Principle [9]  that $\widehat
C=K\cup\widehat K'$.  Now if $q\in V\setminus\widehat{K'}_{\cO(V)}$,
there exists a $\varphi\in\cO(V)$ such that
$\|\varphi(q)\|>\|\varphi\|_{\widehat{K'}_{\cO(V)}}$. Since we for
arbitrarily large $R\in\RR$  have that $V'=V\cap\overline B_R$ is
polynomially convex, $\varphi$  can be approximated uniformly on
$V'$  by an entire function, and the result follows.
\end{proof}

\label{paths}\begin{lemma} Let $K\subset\CC$  be polynomially
convex, let $p_1,p_2\in K$, and let
$Q=\{q_1,...,q_m\}\subset\CC\setminus K$. Then there exists a
polynomially convex set $K'$  such that $K\subset K'$, such that
$p_1$ and $p_2$  are in the same path-connected component of $K'$,
and such that $Q\subset\CC\setminus K'$.
\end{lemma}
\begin{proof}
Choose $R\in\RR$  such that $K\subset\triangle_R$, and let
$\gamma_i:[0,1]\mapsto\CC\setminus K$  be a path connecting $q_i$
and a point $q\in\CC\setminus\triangle_R$  for $i=1,...,m$.  Let
$\gamma:[0,1]\mapsto\triangle_R$ be a smooth curve such that
$\gamma(0)=p_1$ and such that $\gamma(1)=p_2$, and such that
$\gamma\cap Q=\emptyset$. We assume that $p_1$ and $p_2$ are in
different path-components of $K$, or else the lemma is trivial. Let
$\{t^1_1,t^1_2,...,t^l_1,t^l_2\}$ be be an increasing sequence of
numbers in the closed unit interval such that if we let
$\Gamma=\gamma\setminus\{\gamma((t^k_1,t^k_2))\}_{k=1}^l$, we have
that the following is satisfied:

\

(i) $\gamma(t^k_i)\in K$ for all $k$  and all $i$.\

(ii) $\Gamma\cap\gamma_i=\emptyset$ for all $i=1,...,m$.\

(iii) $\gamma((t^k_1,t^k_2))\cap K=\emptyset$ for all
$k=1,...,l$.

\

This means that $Q$ is not in the polynomial hull of
$K_0=K\cup\Gamma$.  If the polynomial hull of
$K_0\cup\gamma([t^1_1,t^1_2])$  does not intersect $Q$  we define
$K_1=\widehat K_0\cup\gamma([t^1_1,t^1_2])$.  If not we do the
following: The only possibility for $\widehat K_1$  to intersect
$Q$  is for $\gamma(t^1_1)$  and $\gamma(t^1_2)$  to be in same
connected component of $K_0$.  Denote this component $C$.  But
this means that if we let $\mu$  be a path close enough to $C$
that connects $\gamma(t^1_1)$  and $\gamma(t^1_2)$, and define
$K_1=K_0\cup\mu$, then the polynomial hull of $K_1$ will not
intersect $Q$.  Do the same thing for the rest of the intervals,
and the set $K'=\widehat K_l$  will satisfy the claims of the
lemma.
\end{proof}

\section{Disjoint Fatou-Bieberbach Domains Whose Union is Dense in
$\CC^k$} \

In this and in the next section we will let
$A:\CC^k\rightarrow\CC^k$  denote the linear map defined by
$$
A:(z_1,...,z_k)\mapsto(\frac{z_1}{2},...,\frac{z_k}{2}).
$$
Fix a $\r>0$.  By Schwarz Lemma there exists a positive number
$\d(\r)$, and two numbers $r,s\in\RR^+$  with $r^2<s$, such that for
a biholomorphic map $F:\overline B_\r\rightarrow\CC^k$  fixing the
origin we have that
$$
\|F-A\|_{\overline B_\r}<\d(\r)\Rightarrow F\in\Gamma(s,r,\r'),
$$
for some positive $\r'$  smaller than $\r$  (the family
$\Gamma(s,r,\r')$ is defined in Lemma 1). By Theorem 4 then, if
$\{F_j\}\subset Aut_p(\CC^k)$ is a sequence of automorphisms
satisfying
$$
(*) \ \|F_j(z)-A(z-p)-p\|<\d(\r), z\in B_\r(p)
$$
for all $j\in\NN$, then $\Omega_{\{F_j\}}^p$  is biholomorphic to
$\CC^k$.  The notation $\d(\r)$  will be used in the following
proof.

\

\emph{Proof of Theorem 1:}  We will prove the result in the case of
$m=\infty$, and we will indicate at the end of the proof what to do
in the finite case.  Let $\e_j\searrow 0$.  To construct the
domains, we will inductively construct a sequence of attracting
automorphisms. At each step in the construction we will generate one
more basin of attraction, and we will make sure that enough points
gets pulled into all of the basins to ensure the claims of the
theorem. \

Let $p_1=q_1$  be the origin and let $\r_1=\frac{1}{2}$.  We start
our construction by letting $F_1$ be the linear map $A$. \

Having constructed $j$ automorphisms, let the following be the
situation $S_j$: We have constructed automorphisms
$\{F_1,...,F_j\}$, we have chosen two sets of distinct points
$\{p_1,...,p_j\}$  and $\{q_1,...,q_j\}$, and a set of positive
numbers $\{\r_1,...,\r_j\}$.  For each $\r_i$  there is a
corresponding $\d(\r_i)$.  The following are satisfied: \

(a) $\overline B_{\r_i}(q_i)\cap\overline B_{\r_k}(q_k)=\emptyset$
for all $i\neq k$, \

(b) $F(j)(p_i)=q_i$  for $i=1,...,j$,\

(c) $F_i(q_k)=q_k$  for $k\leq i$  for all $i=1,..,j$, \

(d) $\|F_i(z)-A(z-q_k)-q_k\|<\d(\r_k)$  for all $z\in\overline
B_{\r_k}(q_k)$ and $k\leq i$. \

We also assume that $\cup_{i=1}^j\overline B_{\r_i}(q_i)$  is
polynomially convex.  When these points, numbers and automorphisms
are chosen at a certain step, they will stay with us throughout the
construction. Notice that because of $(c)$, $(d)$ and Theorem 4; for
any sequence of automorphisms constructed so as to satisfy these
conditions at each step, we have that the basin of attraction of
each point $q_i$ is a Fatou-Bieberbach domain. \

We will now demonstrate how to construct the automorphism $F_{j+1}$.
In addition to ensuring that the four stated claims are satisfied at
the next step, we must make sure that we get enough points into the
basins to satisfy the other claims of the theorem.

Let $K_j=F(j)^{-1}(\cup_{i=1}^j\overline B_{\r_i}(q_i))$.  Choose
sets of points $T_i=\{t^i_1,...,t^i_{m_i}\}\subset B_{j+1}\setminus
K_j$ for $i=1,...,j+1$  such that the sets $T_i$  are pairwise
disjoint. Make sure that for any $q\in B_{j+1}\setminus K_j^\circ$,
for each $i=1,...,j+1$  there is a $t^i_k$  such that \

(e) $\|t^i_k-q\|<\e_j$. \

Let $\tilde t^i_k=F(j)(t^i_k)$.   Next choose a point
$p_{j+1}\in\CC^k\setminus(\cup_{i=1}^{j+1} T_i\cup K_j)$, and write
$q_{j+1}=F(j)(p_{j+1})$.  We may now choose $\r_{j+1}>0$  such that
the set $\overline B_{\r_{j+1}}(q_{j+1})$  does not contain any of
the $\tilde t^i_k$'s  and does not intersect $F(j)(K_j)$.  Make sure
that $\r_{j+1}$  is small enough so that $F(j)(K_j)\cup \overline
B_{\r_{j+1}}(q_{j+1})$  is polynomially convex.\

For a $\mu>0$  let $B^\mu_i$  denote a $\mu$-neighborhood of
$\overline B_{\r_i}(q_i)$  for $i=1,...,j+1$.  For the construction
of the automorphism we will invoke Theorem 2.3 in [5]. By this
result, if $\mu$  is small enough, for any $\e>0$ there exists an
automorphism $\varphi\in Aut(\CC^k)$ such that the following is
satisfied: \

(f) $\|\varphi(z)-A(z-q_i)-q_i\|<\e$  for $z\in B^\mu_i$ for
$i=1,...,j+1$. \

We may also assume that $\varphi(q_i)=q_i$.  By the same theorem
there exists now an automorphism $\phi\in Aut(\CC^k)$  such that \

(g) $\|\phi(z)-z\|<\e$  for $z\in\overline B_{\r_i}(q_i)$  for
$i=1,...,j+1$, \

(h) $\phi(\tilde t^i_k)\in\varphi^{-1}(B_{\r_i}(q_i))$  for all
$t^i_k\in T_i$  for $i=1,...,j+1$. \

Again we may assume that $\phi(q_i)=q_i$.  If we choose $\e$  small
enough, we see that (a),(b),(c) and (d) is satisfied at the new step
$S_{j+1}$ if we define
$$
F_{j+1}=\varphi\circ\phi.
$$.
Notice that if $\e$  is small enough, we also have \

(i) $F(j+1)(t^i_k)\subset B_{\r_i}(q_i)$  for all $t^i_k\in T_i$ for
$i=1,...,j+1$. \

We have now inductively constructed an infinite sequence of
automorphisms $\{F_j\}$.  As commented on earlier, all the basins
$\Omega_i=\Omega_{\{F_j\}}^{q_i}$ are biholomorphic to $\CC^k$.
Since the basins are clearly disjoint, we have ensured (i).  Let
$\Omega$ denote the union of all the basins, let $q\in
B_j\cap(\CC^k\setminus\Omega)$, and let $\Omega_i$  be an arbitrary
basin.  By (e) and (i) there is a sequence of points
$\{t^i_j\}_{j=1}^\infty\subset\Omega_i$  such that
$\|t^i_j-q\|<\e_j$  for all $j\in\NN$.  This shows that
$q\in\partial\Omega_i$, and we have (ii). \

In the case $m\in\NN$  we prove the theorem in the exact same manner
except that we stop generating new basins at the appropriate step in
the construction. \ $\square$

\

\section{Intersections With Subvarieties}

\

In this section we prove Theorem 2 and Theorem 3.  For both proofs
let $A\colon\CC^k\rightarrow\CC^k$ be the linear map defined in the
previous section. \

\

\emph{Proof of Theorem 2:} \ We will construct a sequence of
automorphisms $\{F_j\}\subset Aut_0(\CC^k)$  such that
$\Omega_{\{F_j\}}^0$  is a Fatou-Bieberbach Domain containing
connected subsets of the affine spaces.  Let $\{p^j_i\}_{i\in\NN}$
be a dense set of points in $L_j$ for all $j\in\NN$, and make the
following induction hypothesis $I_j$: We have automorphisms
$\{F_1,...,F_j\}\subset Aut_0(\CC^k)$, and a set of points
$\{q_1,...,q_j\}$  such that $q_i\in L_i$.  For each $k\leq j$ there
are paths $\gamma^k_{lm}\subset L_k$ connecting $p^k_l$  and $p^k_m$
for $l,m\leq j-k+1$.  The following are satisfied:
\begin{align}
&\mathrm{Each} \ F_i \ \mathrm{is \ a \ composition \ of \ maps \
satisfying} \ (*) \ (\mathrm{section \ 5}), \\
&F(j)(p^k_i)\subset B  \ \mathrm{for} \ i\leq j-k+1, k=1,...,j,\\
&F(j)(\gamma^k_{lm})\subset B \ \mathrm{for} \ l,m\leq j-k+1, k=1,...,j,\\
&F(j)(q_i)\subset\CC^k\setminus\overline B \ \mathrm{for} \
i=1,...,j.
\end{align}
We may assume that $I_1$  is true with $F_1=A$.  Assume now that
$I_j$ is true. We will construct $F_{j+1}$ so as to ensure that we
have $I_{j+1}$. \

We want to make sure that there is a path $\gamma^1_{(j+1)1}\subset
F(j)(L_1)$ connecting the images of $p^1_{j+1}$  and $p^1_1$, while
at the same time we have an automorphism tucking the path into the
unit ball while keeping the $F(j)(q_i)$'s at a distance. Notice that
this will ensure that this is also the case for paths
$\gamma_{(j+1)i}$ and points $p^1_i$  for $i=2,...,j$. Now
$K=F(j)^{-1}(\overline B)$ is a polynomially convex set. Let $l_1$
be the complex line containing $p^1_{j+1}$ and $p^1_1$. If none of
the $q_i$'s lie in $l_1$, it follows from Lemma 2 that we for any
path $\gamma\subset l_1$ connecting $p_{j+1}$  and $p_1$ have that
$$
\{q_1,...,q_j\}\cap\widehat{K\cup\gamma}=\emptyset.
$$
If  $Q=\{q_{i_1},...,q_{i_t}\}\subset l_1$,  Lemma 3 tells us that
there is a polynomially convex compact set $K_1'\subset l_1$
containing $K\cap l_1$, such that $Q\cap K_1'=\emptyset$, and such
that $K_1'$ contains a path $\gamma$  connecting $p^1_{j+1}$ and
$p^1_1$.  Lemma 2 tells us that
$$
\widehat{K\cup K_1'}=K\cup
K_1'\Rightarrow\{q_1,...,q_j\}\cap\widehat{K\cup K_1'}=\emptyset.
$$
Let $K'=F(j)(K\cup K_1')$.  Now choose $s\in\NN$ such that
$A^s(K')\subset B$.  For any $\mu>0$, by Theorem 2.1 in [5] there
exists an automorphism $\sigma\in Aut_0(\CC^k)$ such that
$\|\sigma-id\|_{K'}<\mu$, and such that
$\sigma(F(j)(q_i))\in\CC^k\setminus A^{-s}(B_2)$  for $i=1,..,j$.
Define $\psi_1=A^{s}\circ\sigma$. Make sure that $\mu$
 is small enough for $\psi_1$  to be a composition of maps
 satisfying the condition $(*)$  in the beginning of the previous section.\

Now repeat this procedure for the rest of the indices $i=2,...,j+1$.
That is: Construct automorphisms that tuck the $p^i_{j+2-i}$'s along
with the paths into the unit ball, while keeping the $q_i$'s away.
Call the automorphisms $\psi_i$. We will have $I_{j+1}$  if we
define $F_{j+1}=\psi_{j+1}\circ ...\circ\psi_1$, and choose a
$q_{j+1}\in L_{j+1}$  such that $F(j+1)(q_{j+1})\notin\overline B$.

We have inductively defined a sequence of automorphisms $\{F_j\}$
and we claim that $\Omega=\Omega_{\{F_j\}}^0$  is the
Fatou-Bieberbach Domain that we are after.  It follows from Theorem
4 and the choices of automorphisms that $\Omega$ is biholomorphic to
$\CC^k$.  Let $U_1$  and $U_2$ be connected components of
$\Omega\cap L_j$  for a $j\in\NN$. Since the set
$\{p^j_i\}_{i\in\NN}$  is dense in $L_j$, there is a $p^j_l\in U_1$
and a $p^j_m\in U_2$ and we have a path $\gamma^j_{lm}\subset L_j$
connecting the two points while satisfying
$\gamma^j_{lm}\subset\Omega$.  Thus $U_1=U_2$, and we must have that
$\Omega\cap L_j$  is connected.  Lastly we have that
$q_j\subset\L_j\setminus\Omega$  for all $j\in\NN$. \ $\square$

\

It should now be clear how to the prove Theorem 3.  For each
subvariety we choose a compact exhaustion, and we construct the
sequence of automorphisms such that we for each new step tuck more
and larger compact sets into the ball.  Of course we have to make
sure that the basin is not the whole of $\CC^k$.

\

\emph{Proof of Theorem 3:} \ We will construct a sequence of
automorphisms $\{F_j\}\subset Aut_0(\CC^k)$  such that
$\Omega_{\{F_j\}}^0$  is a Fatou-Bieberbach Domain containing
$V=\cup_{i=1}^\infty V_i$. Let $K^j_i$ be $V_j\cap\overline B_i$,
and let $p\in\CC^k\setminus V\cup\overline B_2$.\

We make the following induction hypothesis $I_j$:  We have a
collection of automorphisms $\{F_1,...,F_j\}$  such that
$F(j)(K^k_{j-k+1})\subset B$ for $k\leq j$, and such that
$F(j)(p)\subset\CC^k\setminus\overline B$. By letting $V_1$  be a
coordinate axis, $I_1$  is satisfied with $F_1=A$. \

Assume that we have $I_j$. $F(j)^{-1}(\overline B)$ is a
polynomially convex compact set, and by repeated use of Lemma 2, the
polynomial hull of $K=F(j)^{-1}(\overline B)\cup K^1_{j+1}\cup
K^2_j\cup ...\cup K^{j+1}_1$ does not contain the point $p$. Let
$s\in\NN$  such that $A^s(F(j)(K))\subset B$. For any $\mu>0$, by
Theorem 2.1 in [5] there is an automorphism $\phi\in Aut_0(\CC^k)$
such that $\|\phi-id\|_{F(j)(K)}<\mu$, and such that
$\phi(F(j)(p))\in\CC^k\setminus A^{-s}(B_2)$. Now let
$F_{j+1}=A^s\circ\phi$.  This gives us $I_{j+1}$. Make sure that
$\mu$  is chosen such that $F_{j+1}$  is a composition of maps
satisfying the condition $(*)$  in the beginning of the previous
section. \

We have inductively constructed a sequence of automorphisms
$\{F_j\}$.  It follows from Theorem 4 that
$\Omega=\Omega_{\{F_j\}}^0$ is biholomorphic to $\CC^k$. It is clear
from the construction that all the $V_j$'s will be in the basin, and
also that the point $p$ will not be.  Thus $\Omega$ is the desired
Fatou-Bieberbach Domain. \ $\square$

\

\section{Attracting Basins}

A natural question is the following: Are all Fatou-Bieberbach
domains attracting
 basins for sequences of automorphisms?  We can not answer this in
 general, but in the case of Fatou-Bieberbach domains that are
 also Runge we have the following results:

\

\begin{proposition}
There exists a Fatou-Bieberbach domain $\Omega$  such that there is
no $F\in Aut(\CC^k)$  with
 $\Omega=\Omega_{\{F\}}^p$  for a $p\in\CC^k$.
\end{proposition}
\begin{proof}
For any polynomially convex compact set $K\subset\CC^k$, there
exists a Fatou-Bieberbach domain
 $\Omega$  lying dense in $\CC^k\setminus K$  (See [10] for the strictly convex case).  Now, let
 $U\subset\subset\CC^k$  be an open set such that $(\overline U)^\circ= U$  with $\overline U$  polynomially convex, and let $\Omega$  be a Fatou-Bieberbach
 domain that lies dense in the complement of $\overline U$.  If $\Omega$  is to be
 the attracting basin for an automorphism of $\CC^k$, it is clear that this
 automorphism will have to be an automorphism of $U$.  So if we choose
 $U$  such that no automorphism of $U$  extends holomorphically to an automorphism
 that is attracting at a point outside of $\overline U$, we know that the Fatou-Bieberbach domain in question cannot
 be the attracting basin of an automorphism of $\CC^k$.  And since no automorphism
 of the unit ball extends holomorphically to such an attracting fix-point automorphism,
 we can let $U$  be the unit ball.
\end{proof}
To prove that all Runge Fatou-Bieberbach domains are
sequence-attracting basins we will need the following Lemma: \

\begin{lemma}
Let $\Omega$  be a Runge Fatou-Bieberbach domain.  For a compact
set $K\subset\Omega$,
 a bounded open set $U$  such that $K\subset U$, and an $\e>0$, there exists a
 $\Phi\in Aut(\CC^k)$  such that
$$
\|\Phi(z)-z\|_{K}<\e\quad\Phi(\CC^k\setminus\Omega)\cap U=\emptyset.
$$
\end{lemma}

\begin{proof}
Let $\Psi\colon\CC^k\rightarrow\Omega$  be a Fatou-Bieberbach map.
By [1] there
 exists a sequence of automorphisms $\{F_j\}$  such that $F_j\rightarrow\Psi$  uniformly
 on compacts in $\CC^k$.  So $F_j\circ\Psi^{-1}\rightarrow id$  on $K$.  This
 means that we can let $\tilde\Psi$  be a Fatou-Bieberbach map $\tilde\Psi\colon\Omega\rightarrow\CC^k$
  such that $\|\tilde\Psi(z)-z\|<\frac{\e}{2}$  for all $z\in K$.  By [1] and Corollary 5.3 in [4, p.141]  there
 exists a sequence of automorphisms $\{\Phi_j\}$  such that $\Phi_j\rightarrow\tilde\Psi$
  uniformly on compacts in $\Omega$. So for a large enough $j$  we have the $\e$-estimate.
 And by Corollary 5.3 in [4, p.141] we have that $\|\Phi_j\|\rightarrow\infty$  uniformly
 on $\CC^k\setminus\Omega$, so the result follows with $\Phi=\Phi_m$  for a large enough $m$.
\end{proof}
\

\begin{proposition}
Let $\Omega$  be a Runge Fatou-Bieberbach domain.  For any
$p\in\Omega$  there
 exists a sequence $\{\varphi_j\}\subset Aut_p(\CC^k)$  such that $\Omega=\Omega_{\{\varphi_j\}}^p$.
Moreover, in the terminology from Lemma 1, we may assume that we
have $[\varphi_i(z+p) - p]\in\Gamma(s,r,\r)$  with $r^2<s$  for all
$i\in\NN$.
\end{proposition}
\begin{proof}
We may assume $p=0$, and we let $r_j\searrow 0$  such that
$\overline B_{r_0}\subset\Omega$.  Choose a compact exhaustion
$$
K_0\subset K_1\subset\cdot\cdot\cdot\subset
K_j\subset\cdot\cdot\cdot
$$
of $\Omega$  where $K_0=\overline B_{r_0}$. Let $F_0=A$, where $A$
is the linear map defined in Section 5. \

Now, make the following induction hypothesis $I_j$: We have
automorphisms $\{F_0,...,F_j\}\subset Aut_0(\CC^k)$ such that the
following are satisfied \

(a) $F(j)(K_j)\subset B_{r_ j}$,
\

(b) $F(j)(\CC^k\setminus\Omega)\cap\overline B_{r_0}=\emptyset$. \

$I_0$ is obviously true if $r_0$  is chosen to be small enough. \

Let $r\geq j+1$  such that $B_{r_0}\subset F(j)(K_r)$.  There exists
an $s\in\NN$ such that
$$
A^s(F(j)(K_r))\subset\subset B_{r_{j+1}}.
$$
Let $U$  be a bounded open set such that
$A^{-s}(B_{r_0})\subset\subset U$. Since we have that
$F(j)(K_r)\subset F(j)(\Omega)$  which is a Fatou-Bieberbach domain,
Lemma 4 gives us a $\phi_j\in Aut(\CC^k)$ such that \

(c) $\phi_j\approx id$  on $F(j)(K_r)$, \

(d) $\phi_j(\CC^k\setminus(F(j)(\Omega))\cap U=\emptyset$. \

We may also assume that $\phi_j(0)=0$.  Then we can define
$F_{j+1}=A^s\circ\phi_j$, and $I_{j+1}$  follows.  It is now clear
that $lim_{j\rightarrow\infty} F(j)(z)\rightarrow p$  uniformly on
compacts in $\Omega$, and it is clear that for any
$z\in\CC^k\setminus\Omega$, we do not have convergence of $F(j)(z)$.
Lastly, since $A\circ\phi$  can be made arbitrary close to $A$  on
$\overline B_{r_0}$  we may assume that $F_{j+1}$  is a composition
of maps $\varphi_i$  all elements in $\Gamma(s,r,r_0)$.
\end{proof}

\begin{proposition} Let $\{\Omega_j\}$  be an increasing sequence
of Fatou-Bieberbach domains in $\CC^k$  that are all Runge.  Then
$\Omega=\cup_{j=1}^\infty\Omega_j$  is biholomorphicaly equivalent
to $\CC^k$.
\end{proposition}
\begin{proof}
Let $\{K_j\}$  be an increasing sequence of compact sets that
exhausts $\Omega$  such that $K_j\subset\Omega_j$  for all
$j\in\NN$,  and let $B_j$  denote the ball with radius $j$  in
$\CC^k$.  We will inductively construct a sequence of
biholomorphisms that converges to a biholomorphism
$\Phi\colon\Omega\rightarrow\CC^k$. To start the induction let
$\varphi_1$ be any biholomorphism
$\varphi_1\colon\Omega_1\rightarrow\CC^k$, and assume that we have
constructed maps $\varphi_j\colon\Omega_j\rightarrow\CC^k$   for
$j=1,2,..,k$. What we want to do is to construct $\varphi_{k+1}$
such that it is very close to $\varphi_k$  on $K_k$, and such that
their inverses are very close on $\overline B_k$.  For any $\e$ we
can, by the same argument as in the proof of Lemma 4, assume that we
have a biholomorphism $\phi_{k+1}\colon\Omega_{k+1}\rightarrow\CC^k$
such that \

(a) $\|\phi_{k+1}(z)-z\|<\e$  for all $z\in K_k$,
\

(b) $\|\phi_{k+1}^{-1}(z)-z\|<\e$  for all
$z\in\varphi_k^{-1}(\overline{B_k(h))}$.
\

Here $B_k(h)$  denotes an h-neighborhood of $B_k$.  For any $\d>0$,
by [1] we may assume that we have an $F_{k+1}\in Aut(\CC^k)$ such
that \

(c) $\|F_{k+1}(z)-\varphi_k(z)\|<\d$ for all $z\in K_k(h)$,
\

(d) $\|F_{k+1}^{-1}(z)-\varphi_k^{-1}(z)\|<\d$  for all
$z\in\overline B_k$. \

Now we can define $\varphi_{k+1}=F_{k+1}\circ\phi_{k+1}$, and for
any $\r_{k+1}>0$ we can make sure that \

(e) $\|\varphi_{k+1}(z)-\varphi_k(z)\|<\r_{k+1}, z\in\ K_k$, \

(f) $\|\varphi_{k+1}^{-1}-\varphi_k^{-1}(z)\|<\r_{k+1},
z\in\overline B_k$, \

by letting $\e$  and $\d$  be small enough.  We may now assume that
we have an infinite sequence of biholomorphisms
$\varphi_j\colon\Omega_j\rightarrow\C^n$  that satisfies (e) and (f)
where the sequence $\{\r_j\}$  is chosen to be sumable. In the
terminology from [4] we now have that
$(\varphi_j,\Omega_j)\rightarrow(\Phi,\Omega)$.  And by making sure
that the $\r_j$'s are small enough we can guaranty that $\Phi$  is
not degenerate at every point, which tells us that it is 1-1 onto
its image.  By Theorem 5.2 in [4, p. 140] we also have that
$(\varphi_j^{-1},\CC^k)\rightarrow (\Phi^{-1},\Phi(\Omega))$, so by
the convergence of $\{\varphi_j^{-1}\}$, we must have that
$\Phi(\Omega)=\CC^k$.
\end{proof}

\vskip 5pc

\bf \centerline{Bibliography} \rm

\vskip 2pc

[1] E.Andersen, L.Lempert: {\emph{On the Group of Automorphisms of
$\CC^n$}}, Invent.Math. 110, 371-388, 1992 \

[2] G.T.Buzzard, J.F.Forn\ae ss:{\emph{An Embedding of $\CC$ in
$\CC^2$  with hyperbolic complement}}, Math.Ann.306, 539-546,
1996.

[3] S.Bochner, W.Martin: {\emph{Several Complex Variables}},
Princeton Univ. Press, Princeton, N.J., 1948. \

[4] P.G.Dixon, J.Esterle: {\emph{Michaels Problem and The
Poincare-Fatou-Bieberbach Phenomenon}}, Bull.Amer.Math.Soc. 15,
127-187, 1986 \

[5] F.Forstneric, J-P.Rosay:{\emph{Approximation of Biholomorphic
Mappings by Automorphisms of $\CC^k$}}, Invent.math, 112,323. \

[6] J.E.Forn\ae ss:{\emph{Short $\CC^k$}}, Advanced Studies in Pure
Mathematics, Mathematical Society of Japan, 2003.\

[7] J.E.Forn\ae ss, N.Sibony:{\emph{Some Open Problems in Higher
Dimensional Complex Analysis and Complex Dynamics}, Publicacions
Matematiques 45 (2001), 529-647 \

[8] J.Globevnik:{\emph{On Fatou-Bieberbach Domains}}, Math.Z.229,
91-106, 1998.\

[9] H.Rossi:{\emph{The Local Maximum Modulus Principle}}, Ann.of
Math., 72, 1960\

[10] J-P.Rosay, Walter Rudin:{\emph{Holomorphic Maps from $\CC^n$ to
$\CC^n$}}, American Mathematical Society, Volume 310, Number 1,
November 1988 \

[11] E.F.Wold:{\emph{Proper Holomorphic Embeddings of Finitely and
Infinitely Connected Subsets of $\CC$  Into $\CC^2$}}, Preprint\

\

\centerline{------------------------------------------}

\

Erlend Forn\ae ss Wold\

Ph.d. student, University of Oslo.\

Departement of Mathematics\

University of Michigan\

525 E University\

2074 East Hall\

Ann Arbor, MI 48109-1109\

USA\

e-mail: erlendfw@math.uio.no

\end{document}